\documentclass[12pt,leqno]{article}
\usepackage{amsmath,amssymb,amscd,latexsym,eurosym}
\usepackage[all]{xy}
\newcommand{\C}{{\mathbb{C}}}
\newcommand{\cF}{\overset{\circ}{F}}
\newcommand{\oD}{\overline{D}}
\newcommand{\oF}{\overline{F}}
\newcommand{\oG}{\overline{G}}
\newcommand{\Ha}{\mathbb{H}}
\newcommand{\N}{{\mathbb{N}}}
\newcommand{\Q}{{\mathbb{Q}}}
\newcommand{\R}{{\mathbb{R}}}
\newcommand{\Z}{{\mathbb{Z}}}
\newcommand{\id}{\mathrm{id}}
\newcommand{\ext}{\mathrm{ext}}
\newcommand{\Imm}{\mathrm{Im}\,}
\renewcommand{\iint}{\mathrm{int}\,}
\newcommand{\PSL}{\mathrm{PSL}}
\newcommand{\RRe}{\mathrm{Re}\,}
\newcommand{\tr}{\mathrm{tr}}

\newcommand{\strich}{$'$}
\newcommand{\halb}{\frac{1}{2}}

\newtheorem{definition}{Definition} 
\newtheorem{prop}{Proposition}

\newtheorem{example}{Example} 

\newenvironment{aufzaehl}{\begin{list}{}
{\setlength{\labelwidth}{0.8cm}
\setlength{\leftmargin}{1.1cm}
\setlength{\labelsep}{0.3cm}
\setlength{\rightmargin}{0cm}
\setlength{\itemsep}{0ex}
\setlength{\parskip}{1ex}
\setlength{\topsep}{0ex}
\setlength{\parsep}{0ex}
\setlength{\partopsep}{0ex}
}}
{\end{list}}
\begin{document}
\title{What is a fundamental domain?}
\author{J\"urgen Elstrodt}
\date{}
\maketitle
\begin{abstract}
We show by means of various examples that many of the current definitions of the notion of fundamental domain of a Fuchsian group lack an extra condition ensuring that the domain differs from a measurable fundamental set at most by a null set.\\[0.2cm]
AMS 2000 Mathematics Classification Number: Primary 30F35, secondary 11F06.
\end{abstract}
\section{Introduction}
Throughout the following paper let $\Gamma < \PSL_2 (\R)$ be a Fuchsian group, that is, a discrete subgroup of $\PSL_2 (\R)$, acting on the upper half-plane
\[
\Ha := \{ z \in \C : \Imm z > 0 \}
\]
by linear fractional transformations
\[
z \longmapsto Mz := \frac{az + b}{cz + d}
\]
$(z \in \Ha , M = \pm \left( \begin{smallmatrix} a & b \\ c & d \end{smallmatrix} \right) \in \PSL_2 (\R))$. A {\it fundamental set} for the action of $\Gamma$ on $\Ha$ is usually defined to be a set of representatives of the set of orbits $\Gamma \setminus \Ha := \{ \Gamma z : z \in \Ha \}$. Such a fundamental set in itself will be of little use unless it has some extra structure. Intuitively speaking, a {\it fundamental domain} or {\it fundamental region} of $\Gamma$ is defined to be a ``nice'' subset of $\Ha$ which differs from a fundamental set at most in an inessential way. A look into the textbooks of roughly the last fifty years shows that there is a broad diversity in the definitions of this notion. Most authors demand that a fundamental domain be topologically nice (e.g. open, or closed, respectively) and topologically essentially equal to a fundamental set. However, quite a number of these approaches become insufficient as soon as the (hyperbolic) area measure $\mu \; (d \mu = y^{-2} \, dx \, dy)$ of a fundamental domain is considered. Any two (Lebesgue) measurable fundamental sets of $\Gamma$ have the same area.  The same holds true for so-called measurable fundamental domains which differ from measurable fundamental sets at most by null sets. As is well known, the area of a measurable fundamental domain is an important invariant of $\Gamma$ whenever the area is finite. A critical analysis of many of the common definitions of the notion of ``fundamental domain'' now reveals that a ``fundamental domain'' thus defined may well have an area that is quite different from the area of a measurable fundamental set. This will have undesirable consequences. Recall that the index $[\Gamma : \Delta]$ of a subgroup $\Delta < \Gamma$ is equal to the quotient of the areas of measurable fundamental sets of these groups whenever $\Gamma$ has a measurable fundamental set of finite area. This theorem will be wrong for ``fundamental domains'' whose definition is flawed as indicated above. A similar problem arises with the formation of the Petersson scalar product although it is generally pointed out (and occasionally even ``proved'') that the Petersson product does not depend on the choice of fundamental domain.

Of course, the standard explicit examples of fundamental domains such as the modular triangle, the Poincar\'e normal polygon (also known as the Dirichlet fundamental region), or the Ford fundamental region do not have these deficiencies. The issue of the present note is to point out that many of the common definitions lack an extra condition to the effect that the area of the boundary of a fundamental domain be zero. --- In what follows we construct some unorthodox examples of ``fundamental domains'' satisfying the definitions as given in various sources but which do not have the proper area. According to the broad diversity of the definitions we divide this paper into sections reflecting the case of open, or closed, or arbitrary fundamental domains, respectively.
\section{Open fundamental domains}
Throughout the sequel, $\Gamma$ is a Fuchsian group acting on $\Ha$. Two points $z,w \in \Ha$ are said to be {\it equivalent} under $\Gamma$ if $w = Mz$ for some $M \in \Gamma$. --- Many authors define the notion of fundamental domain as follows.

\begin{definition} \label{def1}
A subset $F \subset \Ha$ is called a {\em fundamental domain} of $\Gamma$ whenever it fulfils the following conditions:
\begin{aufzaehl}
\item[\rm (i)] $F$ is open. 
\item[\rm (ii)] No two distinct points of $F$ are equivalent under $\Gamma$. 
\item[\rm (iii)]  For every $z \in \Ha$ there exists some $M \in \Gamma$ such that $Mz \in \oF$ (closure of $F$ in $\Ha$). 
\end{aufzaehl}
\end{definition}
This definition is given e.g. in \cite{A}, p. 30, \cite{Bea1}, p. 70, \cite{Bu}, p. 19, \cite{Gu}, p. 3, \cite{J}, p. 41, \cite{KM}, p. 2, \cite{Le2}, p. 22. Some writers replace condition (i) by the stronger condition
\begin{aufzaehl}
\item[(i\strich)] $F$ is open and connected;
\end{aufzaehl}
see e.g. \cite{I1}, p. 30, \cite{I2}, p. 38, \cite{Sh}, p. 15, \cite{V}, p. 14. Let us call this variant ``Definition \ref{def1}\strich\,''.

\begin{prop} \label{prop1}
Two fundamental domains of $\Gamma$ in the sense of Definition \ref{def1}\strich\ need not have the same (hyperbolic) area.
\end{prop}

\begin{example} \label{example1}
\rm Let $\Gamma = \PSL_2 (\Z)$ and 
\[ 
D = \{ z \in \Ha : |\RRe z| < \halb \; , \; |z| > 1 \}
\]
the familiar modular triangle; $D$ is a fundamental domain of $\Gamma$ in the sense of Definition \ref{def1}\strich. We now construct quite another ``fundamental domain'': Let $(z_n)_{n \ge 1}$ be an enumeration of all points in $D$ with rational real and imaginary parts, and choose for any $n \in \N$ a small open disc $U_n \subset D$ with centre $z_n$ such that $\sum^{\infty}_{n=1} \mu (U_n) < \varepsilon / 2$, where $0 < \varepsilon < 1$ is chosen arbitrarily small. Connect $z_n$ to the point $2i$ by a (hyperbolic) line segment and enlarge this segment slightly to a small open strip $S_n \subset D$ such that $\sum^{\infty}_{n=1} \mu (S_n) < \varepsilon / 2$. Then $F := \bigcup^{\infty}_{n=1} U_n \cup S_n$ satisfies Definition \ref{def1}\strich, but $\mu (F) < \varepsilon$, whereas $\mu (D) = \frac{\pi}{3}$. --- $F$ even satisfies the following strengthened version (ii\strich) of (ii) which is given in \cite{GGP}, p. 5:
\begin{enumerate}
\item[(ii\strich)] No point of $F$ is equivalent to a different point in $\oF$.
\end{enumerate}
The preceding construction works analogously for any Poincar\'e normal polygon of any Fuchsian group $\Gamma$. Even if $\Gamma$ has a normal polygon of infinite area, our construction yields a ``fundamental domain'' $F$ with $\mu (F) < \varepsilon$ --- clearly a most unsatisfactory state of affairs.

We note two more variants of Definitions \ref{def1}, \ref{def1}\strich: Ford (\cite{Fo}, p. 37) and Lehner (\cite{Le1}, p. 111) leave (i), (ii) unchanged and replace (iii) by the following condition:
\begin{aufzaehl}
\item[(iv)] Every neighbourhood of a boundary point of $F$ contains a point of $\Ha \setminus F$ equivalent to a point of $F$.
\end{aufzaehl}
Behnke and Sommer (\cite{BeS}, p. 477) also maintain (i), (ii) and substitute (v) for (iii):
\begin{aufzaehl}
\item[(v)] $F$ cannot be properly enlarged to an open set $G \supsetneqq F$ satisfying (ii).
\end{aufzaehl}
\end{example}

\begin{prop} \label{prop2}
Even under all the conditions (i)--(v), (ii\,\strich) two fundamental domains of $\Gamma$ may have distinct areas.
\end{prop}

\begin{example} \label{example2}
\rm Recall that a {\it Jordan curve} $\gamma$ is a homeomorphism of a circle to its image, called its {\it trace}, $\tr \gamma \subset \C$. A Jordan curve may have quite unfamiliar properties though it always decomposes $\C \setminus \tr \gamma$ into exactly two components, a bounded one and an unbounded one, called its {\it interior} $\iint \gamma$ and its {\it exterior} $\ext \gamma$, respectively. Moreover, for any Jordan curve $\gamma , \tr \gamma = \partial (\iint \gamma) = \partial (\ext \gamma)$, where $\partial$ means boundary. Though it seems to contradict our imagination, the trace of a Jordan curve may well have positive Lebesgue (area) measure; see \cite{H}, p. 374 f., \cite{K}, \cite{KK}, \cite{Sa}, \cite{ST} for more information on this stunning fact.

To construct the desired example let $\Gamma = \PSL_2 (\Z)$ and $D$ be the open modular triangle as in Example \ref{example1} and choose a Jordan curve $\gamma$ with $\tr \gamma \subset D$ such that $\tr \gamma$ has positive Lebesgue measure. Define $T = \left( \begin{smallmatrix} 1 & 1 \\ 0 & 1 \end{smallmatrix} \right)$, that is, $Tz = z+1 \; (z \in \Ha)$, and put
\[
F := (D \cap \ext \gamma) \cup T (\iint \gamma) \; .
\]
Clearly, $\mu (F) = \mu (D) - \mu (\tr \gamma) < \mu (D)$, and $F$ satisfies (i), (ii), (ii\strich), (iii) since $\tr \gamma = \partial (\ext \gamma) = \partial (\iint \gamma)$. 

To prove (iv) we may right away consider the nontrivial case where $z_0 \in \tr \gamma$ or $z_0 \in T (\tr \gamma)$. Let $z_0 \in \tr \gamma$ and $V$ be a neighbourhood of $z_0$. Since $z_0 \in \partial (\iint \gamma)$, there exists some $z \in V \cap \iint \gamma \subset \Ha \setminus F$, and $z$ is equivalent to $Tz = z + 1 \in F$. The case $z_0 \in T (\tr \gamma)$ is settled similarly. This proves (iv). 

To substantiate (v) assume that the open set $G \supsetneqq F$ satisfies (ii) and that $z_0 \in G \setminus F$. Assume additionally that $z_0$ is an interior point of $G \setminus F$ and let $V \subset G \setminus F$ be an open neighbourhood of $z_0$. By (iii) there exists some $M \in \Gamma$ such that $Mz_0 \in \oF$. Hence $MV \cap \oF \neq \emptyset$ and $MV$ being open it follows that $MV \cap F \neq \emptyset$ since the trace of a Jordan curve cannot have an interior point. But since $V \subset G \setminus F$ this implies that $M \neq \id$ and we obtain a contradiction to condition (ii) for $G$. Hence $G \setminus F$ cannot have an interior point, that is, $F$ is dense in $G$. Hence $F \subsetneqq G \subset \oF$. But this is impossible since every point of $G \setminus F$ must be a boundary point of $F$, and no boundary point of $F$ is an interior point of $\oF$. ---
\end{example}

The construction of Example \ref{example2} works analogously for any
Poincar\'e normal polygon of any Fuchsian group $\Gamma$.

{\bf Open question} Does Proposition \ref{prop2} still hold under the
additional hypothesis that $F$ be connected?

\section{Closed fundamental domains}

There are several distinct definitions of the notion of fundamental domain under the foremost hypothesis that $F$ be closed.

\begin{definition} \label{def2}
A subset $F \subset \Ha$ is called a fundamental domain of $\Gamma$ whenever it satisfies the following conditions:
\begin{aufzaehl}
\item[\rm (i)] $F$ is closed in $\Ha$. 
\item[\rm (ii)] No two distinct points of $\cF$ are equivalent under $\Gamma$. 
\item[\rm (iii)] For every $z \in \Ha$ there exists some $M \in \Gamma$ such that $Mz \in F$.
\end{aufzaehl}
\end{definition}

This definition is given in \cite{KoeK}, first ed., p. 115, second ed. p. 133, \cite{Kob}, p. 100, \cite{N}, p. 220; \cite{Schl}, p. 9 has the same conditions (i), (ii), but omits (iii), which may be assumed to be an oversight.

\begin{prop} \label{prop3}
Two fundamental domains of $\Gamma$ in the sense of Definition \ref{def2} may have distinct areas.
\end{prop}

\begin{example} \label{example3}
\rm Let $\Gamma \neq \{ \id \}$ be a Fuchsian group and $N$ be an arbitrary closed Poincar\'e normal polygon of $\Gamma$. Choose any compact nowhere dense perfect set $P \subset \Ha \setminus N$ of positive Lebesgue measure. (Such a set may be constructed e.g. by a Cantor-like method.) Then $F := N \cup P$ satisfies (i)--(iii) and $\mu (N) < \mu (F)$. ---
\end{example}

Again there are numerous variants of Definition \ref{def2}. Katok (\cite{Ka}, p. 49) demands additionally\footnote{Note that \cite{Ka}, p. 50 expressly states that two fundamental regions of $\Gamma$ have the same area whenever their boundaries are null sets.}:
\begin{aufzaehl}
\item[(iv)] $F = \overline{\cF}$.
\end{aufzaehl}
In the same vein, Miyake (\cite{Mi}, p. 20) demands conditions (i)--(iv), and he also requires that $F$ be a ``connected domain''. It seems not to be clear beyond doubt whether this is intended to mean connectedness of $F$ or of $\cF$. Magnus (\cite{Mag}, p. 57) has a slightly different concept of fundamental domain which we quote as

{\bf Definition 2\strich.} {\it A subset $F \subset \Ha$ is called a fundamental domain of $\Gamma$ whenever there exists an open set $G \subset \Ha$ such that the following conditions hold:
\begin{aufzaehl}
\item[\rm (i\strich)] $F = \oG$ (closure in $\Ha$).
\item[\rm (ii\strich)] No two distinct points of $G$ are equivalent.
\item[\rm (iii\strich)] For every $z \in \Ha$ there exists some $M \in \Gamma$ such that $Mz \in F$.
\end{aufzaehl}}

\begin{prop} \label{prop4}
Under the conditions (i)--(iii) of Definition \ref{def2} plus hypothesis (iv) two fundamental domains of $\Gamma$ may have distinct areas. (Hence the same holds under the conditions of Definition 2\,\strich\ plus hypothesis (iv).)
\end{prop}

\begin{example} \label{example4}
\rm Let $\Gamma = \PSL_2 (\Z)$ and $\oD$ be the closed modular triangle. We go back to Example \ref{example2}, choose a Jordan curve $\gamma$ such that $\tr \gamma \subset \oD$ has positive Lebesgue measure and
\[
\min \{ \RRe z : z \in \tr \gamma \} = - \halb \; , \; \max \{ \RRe z : z \in \tr \gamma \} < \halb \; .
\]
Put
\[
F := (\oD \setminus \iint \gamma) \cup T (\tr \gamma \cup \iint \gamma) \; .
\]
Then $F$ is even connected, satisfies (i)--(iv), and $\mu (F) = \mu (\oD) + \mu (\tr \gamma) > \mu (\oD)$.
\end{example}

{\bf Open Question.} Does Proposition \ref{prop4} hold as well under the additional hypothesis that $\cF$ be connected?

\section{Arbitrary fundamental domains}

We now consider ``fundamental domains'' which are neither assumed to be open nor closed.

\begin{definition} \label{def3}
A subset $F \subset \Ha$ is called a fundamental domain of $\Gamma$ whenever it satisfies the following conditions: 
\begin{aufzaehl}
\item[\rm (i)] Every $z \in \Ha$ has a representative in $F$. 
\item[\rm (ii)] If $z \in F \cap MF$ with $M \in \Gamma, M \neq \id$, then $z \in \partial F$ and $z \in \partial MF$.
\end{aufzaehl}
\end{definition}

This definition is given in \cite{La}, p. 4 and \cite{Maa}, p. 12. By symmetry it is sufficient to demand under (ii) that $z \in \partial F$. Hence Definition \ref{def3} is equivalent to the definition in \cite{Scho}, p. 17. Plainly, Definition \ref{def3} lacks the condition that $F$ be measurable. But even if this hypothesis is added, Example \ref{example4} shows that two ``fundamental domains'' in this sense may have different areas. --- \cite{T}, p. 169 maintains condition (i) of Definition \ref{def3} and substitutes (iii) for (ii):
\begin{aufzaehl}
\item[(iii)] No two distinct points of $\cF$ are equivalent under $\Gamma$.
\end{aufzaehl}
Conditions (i), (iii) admit really strange examples of ``fundamental domains''. E.g., for $\Gamma = \PSL_2 (\Z)$, the domain $\oD$ satisfies (i), (iii), but $F := \oD \cup (\Ha \setminus \Q (i)) , G := \oD \cup (\Ha \cap \Q (i))$ do as well!

Finally, Eichler (\cite{E}, p. 54) defines a fundamental domain of $\Gamma$ to be a fundamental set satisfying in addition the condition
\begin{aufzaehl}
\item[(iv)] $F \subset \overline{\cF}$.
\end{aufzaehl}
As noted above, any two measurable fundamental sets of $\Gamma$ have the same area. But a ``fundamental domain'' in the sense of Eichler's definition need not be Lebesgue measurable as will be shown in Example \ref{example5}.

\begin{example} \label{example5}
\rm Let $\Gamma = \PSL_2 (\Z)$ and 
\[
R := \{ z \in \Ha := - \halb \le \RRe z < \halb , |z| \ge 1 \; , \; \mbox{and} \; |z| > 1 \; \mbox{if} \; \RRe z > 0 \} \; .
\]
$R$ is known to be a fundamental set of $\Gamma$. Now choose a Jordan curve $\gamma$ with $\tr \gamma \subset \overset{\circ}{R}$ such that $\tr \gamma$ has positive Lebesgue measure. As is well known, there exists a nonmeasurable subset $A \subset \tr \gamma$. Define now
\[
F := (R \setminus (A \cup \iint \gamma)) \cup T (A \cup \iint \gamma) \; .
\]
Then $F$ is a fundamental set satisfying (iv), but $F$ is nonmeasurable.
\end{example}

\section{Conclusion}
As we have seen, many of the current definitions of the notion of fundamental domain lack an extra condition ensuring that the domain in question differs from a measurable fundamental set at most by a null set. On the other hand, we quote some sources where such a condition is included. By way of example, H. R\"ohrl in \cite{HC}, p. 587 adds an extra condition to Definition \ref{def2} to the effect that the boundary of $F$ consist of piecewise analytic curves. Similarly, \cite{FK}, p. 206 supplement Definition \ref{def1} by the requirement that the relative boundary of $F$ consist of piecewise analytic arcs that are equivalent in pairs. These definitions clearly exclude our examples based on the existence of Jordan curves whose trace has positive area since already C. Jordan proved that the trace of a rectifiable arc is a null set (see e.g. \cite{BuG} or \cite{Sa}, p. 131 f.). Several other authors just demand that the boundary of $F$ be of Lebesgue area measure zero; see e.g. \cite{Bea2}, p. 204, \cite{Kr}, p. 7, \cite{Mac}, p. M35, \cite{Mas}, p. 28--29.

Mathematisches Institut\\
Westf. Wilhelms-Universit\"at\\
Einsteinstr. 62\\
48149 M\"unster, Germany\\
email: elstrod@uni-muenster.de
\end{document}